\documentclass[11pt,a4paper]{amsart}
\usepackage{amsmath}
\usepackage{amsfonts}

\newtheorem{theorem}{Theorem}

\newtheorem{lemma}{Lemma}

\newtheorem{proposition}{Proposition}

\input{tcilatex}

\title[Polynomials \'{a} la Lehmers and Wilf]{Polynomials \'{a} la Lehmers and Wilf}
\author{Gert Almkvist and Arne Meurman}
\address{Centre for Mathematical Sciences, Mathematics, Lund
  University, Box 118, SE-221 00 Lund, Sweden}
\email{gert@maths.lth.se, arnem@maths.lth.se}

\begin{document}

\maketitle


\section{Introduction.}

In \cite{Alm} a generalized Dedekind sum is defined by%
\begin{equation*}
s(r,h,k)=\frac{k^r}{r+1}\sum_{j=1}^{k-1}B_{r+1}(j/k)((jh/k))
\end{equation*}%
where \ $r$ \ is an even integer, \ $B_{r+1}(x)$ \ the Bernoulli polynomial
and%
\begin{equation*}
((x))=
\begin{cases}
x-[x]-1/2\text{ \ if \ }x\notin \mathbb{Z}\\
0\text{ \ if\ }x\in\mathbb{Z}\\
\end{cases}
\end{equation*}%
Further we define%
\begin{equation*}
A(r,k,n)=\sum_{(h,k)=1}\exp \left\{ \pi is(r,h,k)-2\pi ihn/k\right\} 
\end{equation*}%
We define the generalized Wilf polynomial of degree \ $k$ \ by%
\begin{equation*}
W(r,k,x)=\dprod\limits_{n=1}^{k}(x-A(r,k,n))
\end{equation*}%
T.Dokshitzer \cite{Dok} proved that \ $W(r,k,x)$ \ has integer coefficients if%
\begin{equation*}
(r+1,k)=1\text{ \ and \ }(r+1,\varphi (k))=1
\end{equation*}%
where \ $\varphi $ \ is Euler's totient function. In this note we study \ $%
W(p-3,p,x)$ which has integer coefficients if \ $p$ \ is a prime $\geq 5$.

Let \ $p\geq 5$ \ be a prime and $g$ a primitive root
$\pmod {p^2}$.
Define%
\begin{equation*}
\eta _{n}=\sum_{j=0}^{p-2}\exp \left\{ 2\pi ig^{pj}(g^{p}+pn)/p^2\right\} 
\end{equation*}%
and 
\begin{equation*}
L(p,x)=\dprod\limits_{n=1}^{p}(x-\eta _{n})
\end{equation*}%
the period polynomial \ (see E.Lehmer and D.H.Lehmer \cite{Leh}). We show that \ $%
W(p-3,p,x)$ \ and \ $L(p,x)$ \ agree up to the sign of \ $x.$ In \cite{Leh} the
Lehmers ask if the constant term \ $L(p,0)$ could be even. We computed \ $%
L(1093,0)$ and found it divisible by \ $2^{1102}$ so it is even by some
margin. More precisely we have the following:

Let \ $q$ \ be a prime such that \ $q^{p-1}\equiv 1$ $\pmod{p^{2}}$
(Wieferich condition). Then \ $L(p,x)$ $\pmod{q}$ \ splits into linear
factors.

In a private communication A.Granville noted that \ $W(p-3,p,x)$ \ could be
written as a certain determinant. This is proved here. After reflection
it is the characteristic polynomial of a circulant matrix.

\bigskip

\section{Generalized Wilf polynomials.}

In \cite{Alm} it is proved (using Apostol's reciprocity theorem \cite{Apo}) that if \ $r\leq
p-5$ \ then \ $A(r,p,n)\in \mathbb{Q}[\zeta_p ]$ \ where \ $\zeta_p =\exp
(2\pi i/p).$ This is not the case for \ $r=p-3.$ Instead we have to use \ $%
\omega =\exp (2\pi i/p^{2}).$ Let 
\begin{equation*}
G=\func{Gal}(\mathbb{Q}(\omega )/\mathbb{Q}) 
\end{equation*}%
and \ $g$ \ a primitive root in \ $(\mathbb{Z}/p^{2}\mathbb{Z)}^{\times
}$. Let $\alpha$ be the automorphism of $\mathbb{Q}(\omega)$
determined by $\omega \mapsto\omega^{1+p}$. Then we have
\begin{equation*}
G=U\times V 
\end{equation*}%
where%
\begin{equation*}
U=\left\langle \alpha \right\rangle = \left\{ \omega \mapsto \omega
^{1+jp};j=0,1,...,p-1\right\} 
\end{equation*}%
\begin{equation*}
V=\left\langle \omega \mapsto \omega ^{g^{p}}\right\rangle =\left\{
\omega \mapsto \omega ^{h^{p}};h=1,2,...,p-1\right\} 
\end{equation*}%
Set%
\begin{equation*}
E=\mathbb{Q}(\omega )^{V}=\mathbb{Q}(\rho ) 
\end{equation*}%
where%
\begin{equation*}
\rho =\func{Tr}_{\mathbb{Q}(\omega )/E}(\omega )=\sum_{h=1}^{p-1}\omega ^{h^{p}} 
\end{equation*}%
(it is easy to show that $\rho\notin\mathbb{Q}$ cf. \cite{Leh}).
Then we obtain%
\begin{equation*}
L(p,x)=\dprod\limits_{j=0}^{p-1}(x-\alpha ^{j}(\rho )) 
\end{equation*}%
since $\eta_{mg^p} = \alpha^m(\rho)$.
This has integer coefficients since \ $\rho $ \ is an algebraic integer. We
denote by \ $\mathbb{Z}_{p}$ \ the ring of \ $p$-adic integers. Our main
goal is

\begin{theorem}Let \ $p\geq 5$ \ be a prime. Then%
\begin{equation*}
W(p-3,p,x)=\begin{cases}
L(p,x)\text{ \ if \ }p\equiv \pm 1\text{
  }\pmod{8}\\
 -L(p,-x)\text{ \ if \ }p\equiv \pm 3\text{ }\pmod{8}
\\
\end{cases}
\end{equation*}
\end{theorem}

This will follow from the following two theorems:

\begin{theorem} Let \ $u\in \mathbb{Z}$ \ be such that%
\begin{equation*}
u\equiv \frac{pB_{p-1}}{p-2}\text{ }\pmod{p^{2}\mathbb{Z}_p}
\end{equation*}%
Then%
\begin{equation*}
p^{2}s(p-3,h,p)\equiv uh^{p}\text{ }\pmod{p^{2}} 
\end{equation*}%
for \ $h=1,2,...,p-1$.
\end{theorem}

Note that the von Staudt-Clausen theorem implies that%
\begin{equation*}
u\equiv \frac{p+1}{2}\text{ }\pmod{p} 
\end{equation*}%
so that in particular \ $(u,p)=1.$

\begin{theorem}We have%
\begin{equation*}
p^{2}s(p-3,h,p)\equiv \frac{p^{2}-1}{8}\text{ }\pmod{2} 
\end{equation*}%
for \ $h=1,2,...,p-1.$
\end{theorem}

The proof of Theorem 19 in \cite{Alm} is also valid for \ $r=p-3$ \ and gives
Theorem 3 (one requires \ $r$ \ even and $\ (r+1,p)=1$\;). The proof of
Theorem 16 in \cite{Alm} \ (which uses Apostol's reciprocity theorem
\cite{Apo}) gives%
\begin{equation}
p^{2}s(p-3,h,p)\equiv \frac{pB_{p-1}}{p-1}\left\{ H^{-1}+\frac{1}{p-2}%
H^{p-2}\right\} \text{ }\pmod{p^{2}\mathbb{Z}_{p}} \label{reciprocity}
\end{equation}%
for any integers \ $h,H$ \ such that \ $hH\equiv 1$ $\pmod{p}$.

\begin{lemma}Let \ $a\in \mathbb{Z}_{p}$ such that \ $-2a\equiv 1$ $%
\pmod{p\mathbb{Z}_{p}}$. Define for $H\in\mathbb{Z}$, $H\not\equiv 0 \pmod{p}$,
\begin{equation*}
f(H)=H^{-1}+aH^{p-2}\text{ }\pmod{p^{2}\mathbb{Z}_{p}} 
\end{equation*}%
Then 
\begin{equation*}
f(H+bp)\equiv f(H)\text{ }\pmod{p^{2}\mathbb{Z}_{p}} 
\end{equation*}
for all $b\in\mathbb{Z}$.
\end{lemma}

\textbf{Proof. }It is enough to show 
\begin{equation*}
f(H(1+p))\equiv f(H)\text{ }\pmod{p^{2}\mathbb{Z}_{p}} 
\end{equation*}%
since the general statement follows from this by iteration. We have 
\begin{align*}
f(H(1+p))&-f(H)=\frac{1}{H(1+p)}+aH^{p-2}(1+p)^{p-2}-\frac{1}{H}-aH^{p-2} \\
&=\frac{1}{H(1+p)}\left\{ 1+aH^{p-1}(1+p)^{p-1}-1-p-a(1+p)H^{p-1}\right\} \\
&\equiv  \frac{1}{H(1+p)}\left\{ -p+aH^{p-1}(1-p-1-p)\right\} \\
&\equiv  \frac{p}{H(1+p)}\left\{ -1-2aH^{p-1}\right\} \\
&\equiv \frac{p}{H(1+p)}%
\left\{ H^{p-1}-1\right\} \equiv 0\text{ }\pmod{p^{2}\mathbb{Z}_{p}}  
\end{align*}
\hfill $\Box$

\medskip
\textbf{Proof of Theorem 2. }\ As \ $H^{p}\equiv H$ \ $\pmod{p}$, \
using Lemma 1, we may replace \ $H$ \ by \ $H^{p}$ \ in the right hand side
of \eqref{reciprocity}.
Hence%
\begin{align*}
p^{2}s(p-3,h,p)&\equiv \frac{pB_{p-1}}{p-1}\left\{\frac{1}{H^{p}}+\frac{1}{p-2}%
H^{p(p-2)}\right\} \\
&\equiv \frac{pB_{p-1}}{p-1}\frac{1}{H^{p}}\left\{ 1+\frac{1}{p-2}%
H^{p(p-1)}\right\} \\
&\equiv \frac{pB_{p-1}}{p-1}\frac{p-1}{p-2}\frac{1}{H^{p}}%
\equiv uh^{p}\text{ }\pmod{p^{2}\mathbb{Z}_{p}} 
\end{align*}%
Here we used the fact that if \ $hH\equiv 1$ $\pmod{p}$ \ then \ $
h^{p}H^{p}\equiv 1$ $\pmod{p^{2}}$.

\hfill $\Box$

\bigskip
With Theorem 2 and 3 we can now evaluate \ $A(p-3,p,n).$ Note that%
\begin{equation*}
\exp (\pi i/p^{2})=-\omega ^{(p^{2}+1)/2} 
\end{equation*}%
Thus%
\begin{align*}
A(p-3,p,n)&=\sum_{h=1}^{p-1}\exp \left\{ \frac{\pi i}{p^{2}}%
p^{2}s(p-3,h,p)-2\pi i\frac{hn}{p}\right\} \\
&=\sum_{h=1}^{p-1}(-1)^{p^{2}s(p-3,h,p)}\omega ^{(p^{2}+1)/2\cdot uh^{p}-phn} \\
&=(-1)^{(p^{2}-1)/8}\sum_{h=1}^{p-1}\omega ^{(v-pn)h^{p}} \\
&=(-1)^{(p^{2}-1)/8}\func{Tr}_{\mathbb{Q}(\omega )/E}(\omega ^{v-pn}) 
\end{align*}%
where%
\begin{equation*}
v\equiv \frac{p^{2}+1}{2}u\text{ }\pmod{p^{2}} 
\end{equation*}

With \ $\alpha $ \ the automorphism of \ $\mathbb{Q}(\omega )$ \ determined
by \ $\omega \mapsto \omega ^{1+p}$ \ we now obtain the action of \ $%
U\simeq \func{Gal}(E/\mathbb{Q})$ on the \ $A(p-3,p,n)$%
\begin{align*}
\alpha (A(p-3,p,n))&=(-1)^{(p^{2}-1)/8}\func{Tr}_{\mathbb{Q}(\omega )/E}(\omega
^{(1+p)(v-pn)})\\ 
&=(-1)^{(p^{2}-1)/8}\func{Tr}_{\mathbb{Q}(\omega )/E}(\omega
^{v-p(n-v)})\\
&=A(p-3,p,n-v) \label{galoisaction}
\end{align*}
Thus the $\left\{ A(p-3,p,n);n=0,1,...,p-1\right\} $ \ form one orbit under
the action of \ $U$. If \ $n$ \ is such that%
\begin{equation*}
v-pn\equiv v^{p}\text{ }\pmod{p^{2}} 
\end{equation*}%
then%
\begin{align*}
A(p-3,p,n)&=(-1)^{(p^{2}-1)/8}\func{Tr}_{\mathbb{Q}(\omega )/E}(\omega
^{v^{p}})\\
&=(-1)^{(p^{2}-1)/8}\func{Tr}_{\mathbb{Q}(\omega )/E}(\omega )=\pm \rho 
\end{align*}
and Theorem 1 follows.

\begin{theorem}Let \ $q$ \ and \ $p$ \ be primes satisfying the
Wieferich condition%
\begin{equation*}
q^{p-1}\equiv 1\text{ }\pmod{p^{2}} 
\end{equation*}%
Then \ $L(p,x)$ $\pmod{q}$ \ splits into linear factors.
\end{theorem}

\textbf{Proof.} Let $R$ be the integral closure of $\mathbb{Z}$ in
$E$, the ring of algebraic integers in $E$. Choose an element $\gamma$
of (multiplicative) order $p^2$ in the finite field
$\mathbb{F}_{q^{p-1}}$. Let
\begin{equation*}
\phi:\mathbb{Z}[\omega]\to \mathbb{F}_{q^{p-1}}
  \end{equation*}
  denote the homomorphism determined by $\phi(\omega) = \gamma$.
Then $R\subset \mathbb{Z}[\omega]$, since
  $\mathbb{Z}[\omega]$ is the ring of algebraic integers in
  $\mathbb{Q}(\omega)$, see e.g. \cite{L2} Chapter IV, Theorem 3.
  We shall show
  that
  $\phi(R) = \mathbb{F}_q$. Certainly $\phi(R)$ is an intermediate
  field
  \begin{equation*}
\mathbb{F}_q \subset \phi(R) \subset \mathbb{F}_{q^{p-1}},
    \end{equation*}
so $\func{Gal}(\phi(R)/\mathbb{F}_q)$ is a cyclic group of order
dividing $p-1$. By \cite{L1} Chapter VII, Proposition 2.5
$\func{Gal}(\phi(R)/\mathbb{F}_q)$ is a subquotient of
$\func{Gal}(E/\mathbb{Q})$, the quotient of the decomposition group by
the inertia group. As $\func{Gal}(E/\mathbb{Q})$ is cyclic of order
$p$, the only possibility is that $\func{Gal}(\phi(R)/\mathbb{F}_q) =
1$ so that $\phi(R) = \mathbb{F}_q $.

As $\eta_n\in R$ we have $\phi(\eta_n)\in\mathbb{F}_q$ so that
\begin{equation*}
\phi(L(p,x)) = \prod_{n=1}^p (x-\phi(\eta_n))
\end{equation*}
and the theorem is proved.

\hfill $\Box$

\medskip
\noindent
\textbf{Remark }If \ $q,p$ \ is a Wieferich pair Theorem 4 implies that \
often $q$ \ divides the constant term \ $L(p,0)$ \ to some
high power. We give a small table of this phenomenon%
\begin{equation*}
\begin{tabular}{|l|l|}
\hline
$p$ & factor of \ $L(p,0)$ \\ \hline
$11$ & $3^{5}$ \\ \hline
$13$ & $23^{3}$ \\ \hline
$43$ & $19^{4}$ \\ \hline
$47$ & $53^{2}$ \\ \hline
$59$ & $53^{2}$ \\ \hline
$71$ & $11^{4}$ \\ \hline
$79$ & $31^{5}$ \\ \hline
$97$ & $107^{4}$ \\ \hline
$103$ & $43^{4}$ \\ \hline
$113$ & $373^{4}$ \\ \hline
$137$ & $19^{14}$ \\ \hline
$331$ & $71^{7}$ \\ \hline
$863$ & $13^{80}$ \\ \hline
$1093$ & $2^{1102}$ \\ \hline
\end{tabular}%
\end{equation*}%

\vskip 7mm
\section{Granville's determinant.}

In this section all matrix indices shall be interpreted in
$\mathbb{Z}/p\mathbb{Z}$ with representatives $\{1,2,\dots,p\}$.

\begin{theorem}Let \ $A=(a_{m,n})$ \ be the \ $p\times p$-matrix
with%
\begin{equation*}
a_{m,n}=
\begin{cases}
-x\text{ \ if \ }m+n\equiv 0 \pmod{p}\\
\exp \left\{ 2\pi i(m+n)^{p}/p^{2}\right\}\text{\ otherwise \ }\\
\end{cases}
\end{equation*}%
Then%
\begin{equation*}
\det (A)=(-1)^{(p+1)/2}L(p,x)
\end{equation*}
\end{theorem}

\noindent
\textbf{Example: }Let \ $p=7.$ Then with \ $\omega =\exp (2\pi i/49)$%
\begin{equation*}
A=%
\begin{vmatrix}
\omega ^{30} & \omega ^{31} & \omega ^{18} & \omega ^{19} & \omega ^{48} & -x
& \omega  \\ 
\omega ^{31} & \omega ^{18} & \omega ^{19} & \omega ^{48} & -x & \omega  & 
\omega ^{30} \\ 
\omega ^{18} & \omega ^{19} & \omega ^{48} & -x & \omega  & \omega ^{30} & 
\omega ^{31} \\ 
\omega ^{19} & \omega ^{48} & -x & \omega  & \omega ^{30} & \omega ^{31} & 
\omega ^{18} \\ 
\omega ^{48} & -x & \omega  & \omega ^{30} & \omega ^{31} & \omega ^{18} & 
\omega ^{19} \\ 
-x & \omega  & \omega ^{30} & \omega ^{31} & \omega ^{18} & \omega ^{19} & 
\omega ^{48} \\ 
\omega  & \omega ^{30} & \omega ^{31} & \omega ^{18} & \omega ^{19} & \omega
^{48} & -x%
\end{vmatrix}%
\end{equation*}%
We get%
\begin{equation*}
\det (A)=x^{7}-21x^{5}-21x^{4}+91x^{3}+112x^{2}-84x-97=L(7,x)
\end{equation*}

\textbf{Proof.} \ Let \ $E_{i,j}$ be the $\ p\times p-$matrix with \ $1$ \
in place \ $(i,j)$ \ and \ $0$ \ otherwise. Let \ 
\begin{equation*}
B=\sum_{i\neq j}\omega ^{(i-j)^{p}}E_{i,j}
\end{equation*}%
and%
\begin{equation*}
F=\sum_{j=1}^{p}E_{j,p-j}
\end{equation*}%
Let \ $C=xI-B$. \ Then we have%
\begin{equation*}
C=-FA
\end{equation*}%
Since \ $\det (F)=(-1)^{(p-1)/2}$ \ and the dimension \ $p$ \ is odd we get%
\begin{equation*}
\det (C)=(-1)^{(p+1)/2}\det (A)
\end{equation*}%
so \ $\det (A)$ \ is equal to the characteristic polynomial of the
circulant\ $B$ up to sign. We find the eigenvalues of \ $B$.

\begin{proposition} (cf. \cite{Ait}, page 123) For \ $k=0,1,...,p-1$ \ the vector%
\begin{equation*}
v_{k}=\sum_{j=1}^{p}\omega ^{pjk}e_{j}
\end{equation*}%
where \ $e_{j}$ \ is the column vector with \ $1$ \ at place \ $j$ \ and \ $0
$ \ elsewhere, is an eigenvector of \ $B$ \ with eigenvalue%
\begin{equation*}
\rho _{k}=\sum_{j=1}^{p-1}\omega ^{(1-pk)j^{p}}
\end{equation*}
\end{proposition}

\textbf{Proof. }Introduce the matrix%
\begin{equation*}
T=\sum_{j=1}^{p}E_{j+1,j}
\end{equation*}%
Then%
\begin{equation*}
B=\sum_{j=1}^{p-1}\omega ^{j^{p}}T^{j}
\end{equation*}%
and hence each eigenvector of \ $T$ \ is an eigenvector of \ $B.$ Now%
\begin{equation*}
Tv_{k}=\omega ^{-pk}v_{k}
\end{equation*}%
so%
\begin{align*}
Bv_{k}&=\sum_{j=1}^{p-1}\omega ^{j^{p}}T^{j}v_{k}=\sum_{j=1}^{p-1}\omega
^{j^{p}}\omega ^{-pjk}v_{k} \\
&=\sum_{j=1}^{p-1}\omega ^{(1-pk)j^{p}}v_{k}=\rho _{k}v_{k} 
\end{align*}%
since \ $j^{p}\equiv j$ $\pmod{p}$, completing the proof of the
Proposition.

As the $\rho_k$ form one orbit of $U$ acting on $\rho = \rho_0$,
\begin{equation*}
\det (C)=\det (xI-B)=\dprod\limits_{k=0}^{p-1}(x-\rho _{k})=L(p,x)
\end{equation*}%
finishes the proof of Theorem 5.

\end{document}